\documentclass[letterpaper,11pt]{article}
\usepackage{amsmath, amsthm, amssymb}    
\usepackage{bridges}			
\usepackage{graphicx}			
\usepackage{tikz}
\usetikzlibrary{positioning}
\usetikzlibrary { decorations.pathmorphing, decorations.pathreplacing, decorations.shapes, }
\usetikzlibrary{arrows.meta}

\usepackage[colorlinks=true, urlcolor=blue, citecolor=black, linkcolor=black]{hyperref}  
\usepackage{subcaption}			

\newcommand{\1}{\times}
\newcommand{\0}{\mathop{\cdot}}
\newcommand{\da}{\downarrow}

\newcommand{\R}{\mathbb{R}}

\newcommand{\floor}[1]{\left\lfloor #1\right\rfloor}

\definecolor{mycolor}{RGB}{252,150,255}

\newtheorem{theorem}{Theorem}

\theoremstyle{definition}
\newtheorem{definition}{Definition}

\newtheorem{example}{Example}

\title{The Music and Mathematics of Maximal Evenness in Graphs}

\author{Neal Bushaw,\textsuperscript{1} Brent Cody,\textsuperscript{2}  Luke Freeman\textsuperscript{3}, and Tobias Whitaker\textsuperscript{4}
\vspace{10pt}\\
\textsuperscript{1}Virginia Commonwealth University, Richmond, Virginia, USA; nobushaw@vcu.edu\\
\textsuperscript{2}Virginia Commonwealth University, Richmond, Virginia, USA; bmcody@vcu.edu\\
\textsuperscript{3}Virginia Commonwealth University, Richmond, Virginia, USA; freemanln@vcu.edu\\
\textsuperscript{4}University of Richmond, Richmond, Virginia, USA; twhitak2@richmond.edu \\
} 

\date{}					

\begin{document}

\maketitle

\thispagestyle{empty}

\begin{abstract}

We use the concept of electric potential energy from physics, the mathematical field of graph theory, and the notion of majorization to study maximal evenness in a broader mathematical context than what was previously possible, so that we can go beyond the well-known one-dimensional maximally even sets into higher dimensional and more geometrically complex territory. We investigate musical connections between certain generalizations of maximally even sets, one of the oldest Puerto Rican musical traditions of African origin called bomba, and with certain scales ranging from the familiar to the esoteric.
\end{abstract}


\section*{Introduction}



\begin{figure}[h!tbp]
\centering
\begin{minipage}[b]{0.15\textwidth} 
	\includegraphics[width=\textwidth]{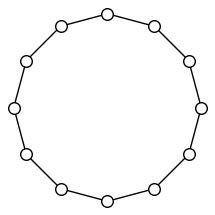}
        	\subcaption{} 
        	\label{figure_12cycle}
\end{minipage}
\begin{minipage}[b]{0.15\textwidth} 
	\includegraphics[width=\textwidth]{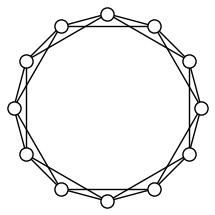}
        	\subcaption{} 
        	\label{figure_cyclesquared}
\end{minipage}
\begin{minipage}[b]{0.45\textwidth} 

\ \ \ \ \begin{tikzpicture}
\tikzset{highlight/.style={preaction={
draw,mycolor,-,
double=mycolor,
double distance=10\pgflinewidth,opacity=0.3
}}}

\foreach \i [count=\xi from 0, remember=\i as \iprev] in {0,...,7}{
  \node[draw,circle,inner sep=0, minimum size=1.5mm] (\i) at ({1.5*cos(360*\i/8-360/7)},{0.75*sin(360*\i/8-360/7)}) {};

  \ifthenelse{\i>0}{
    \draw[-,highlight] (\iprev) -- (\i);
  }{}
}

\foreach \i [count=\xi from 0, remember=\i as \iprev] in {0,...,7}{
  \node[draw,circle,inner sep=0, minimum size=1.5mm, fill=white] (a\i) at ({1.5*cos(360*\i/8-360/7)},{0.75*sin(360*\i/8-360/7)+0.7}) {};

  \draw[-] (\i) -- (a\i);

  \ifthenelse{\i>0}{
    \draw[-,highlight] (a\iprev) -- (a\i);
  }{}

}

\draw[-,highlight] (0) -- (a7);
\draw[-,highlight] (a0) -- (7);


\node[] () at (2,0.3) {$\cong$};

\foreach \i [count=\xi from 0, remember=\i as \iprev] in {0,...,15}{
  \node[draw,circle,inner sep=0, minimum size=1.5mm] (\i) at ({1.2*cos(360*\i/16)+3.8},{1.2*sin(360*\i/16)+0.3}) {};

  \ifthenelse{\i>0}{
    \draw[-,highlight] (\iprev) -- (\i);
  }{}

  \ifthenelse{\i=15}{
    \draw[-,highlight] (0) -- (\i);
  }{}

}

\foreach \i [evaluate={\opp=int(\i+8);}] in {0,...,7}{
  \draw[-] (\i) -- (\opp);
}
\end{tikzpicture}

        	\subcaption{} 
        	\label{figure_mobiusladder}
\end{minipage}
\begin{minipage}[b]{0.15\textwidth} 
	\includegraphics[width=\textwidth]{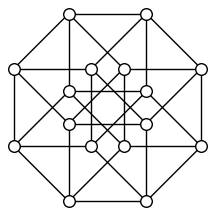}
        	\subcaption{} 
        	\label{figure_hypercube}
\end{minipage}
~ 
\caption{Various finite graphs:  (a) a $12$-cycle $C_{12}$, (b) a $12$-cycle squared $C_{12}^2$, (c) two isomorphic representations of the M\"obius ladder $M_{16}$, and (e) a hypercube.  }
\label{figure_graphs}
\end{figure}
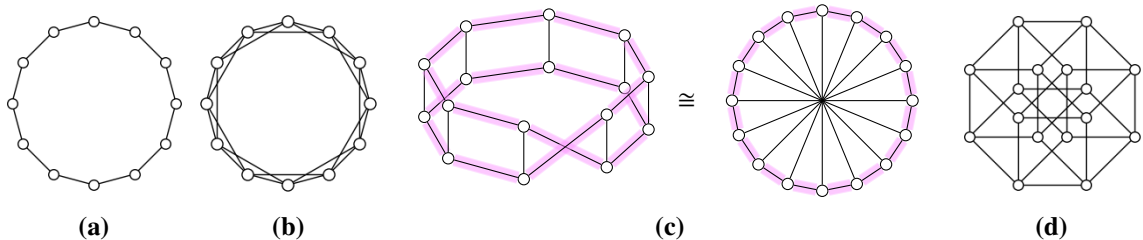

Maximally even sets, which arose as part of Clough and Douthett's study \cite{CloughDouthett} of musical scales, and which have gained in popularity more recently due to the rhythmic connections discovered by Toussaint \cite{Toussaint04}, manifest in musical traditions across cultures in the form of interesting scales and rhythms.  In a remarkable turn of events, Clough and Douthett's research in music theory has found applications in mathematics \cite{MR2408358,MR2212108}, computer science \cite{MR2069421}, mathematical physics \cite{MR1401228}, and even in the design of particle accelerators \cite{Bjorklund2004TheTO}. Up until now, maximally even sets have always been assumed to exist within a ``cyclic'' context. For example, scales are often viewed as particular subsets of the twelve standard pitch classes arranged around a circle, and repeating rhythms are often visualized as a collection of beats and rests arranged in a circular pattern. It is natural to view the cyclic contexts of pitch classes and rhythms as being graphs (i.e. collections of vertices and edges): the twelve pitch classes correspond to the vertices of the $12$-cycle $C_{12}$ (Figure \ref{figure_12cycle}) and a rhythm of length twelve can be viewed as a particular subset of $C_{12}$. In this article we address the following two questions. First, to what extent can the definitions and techniques used to study maximal evenness in the cyclic context be generalized to other graphs, such as squares of cycles (Figure \ref{figure_cyclesquared}), M\"obius ladders (two representations of which are given in Figure \ref{figure_mobiusladder}), or hypercubes (Figure \ref{figure_hypercube})? Second, are there musical connections to be found by considering maximal evenness in graphs other than cycles?

\section*{Maximal Evenness}


Suppose we want to write a repeating rhythm for a single drum which contains five onsets (i.e. hits on the drum) and seven rests, for a total of twelve beats in our measure, where each onset is denoted by ``$\1$'' and each rest by ``$\cdot$.'' Further suppose that we want the onsets to be spaced out ``as evenly as possible'' within the twelve available beats. As a first attempt, we might begin by calculating $\frac{12}{5}=2.4$, and since $2.4$ is closer to $2$ than to $3$, we might try placing an onset on every two beats like this: $[\1\0\1\0\1\0\1\0\1\0\0\0]$. But one quickly realizes that the onsets can be ``more evenly spaced'' by writing $[\1\0\1\0\1\0\0\1\0\1\0\0]$, which is indeed a ``maximally even'' configuration of onsets. Since this is a repeating rhythm, it can be visualized as a particular set of vertices $A=\{0,2,4,7,9\}$ in the graph $C_{12}$ (see Figure \ref{figure_equilibrium}), and $A$ is called a \emph{maximally even set} or a \emph{Euclidean rhythm} \cite{MR2512671, Toussaint04}.

There are many equivalent algorithms that generate maximally even sets \cite{CloughDouthett, MR2512671, MR2408358}, but perhaps one of the most straightforward, which we take as our definition, is that involving $J$-representations, which is due to Clough and Douthett \cite{CloughDouthett}. Before stating the definition, let us consider an example of a $J$-representation of the maximally even set $A=\{0,2,4,7,9\}$ in $C_{12}$, which one can think of as being a formula for $A$. Clough and Douthett defined the notation $J^0_{12,5}=\left\{\floor{\frac{12i}{5}}:i=0,1,2,3,4\right\}$,
and with a little simplification we see that $A=J^0_{12,5}$. Let us note that one can view the definition of $J^0_{12,5}$ as being a ``fix'' of the inadequate first attempt made above in which we computed $\frac{12}{5}=2.4\sim 2$ and placed onsets every two beats.

\begin{definition}\label{definition_maximally_even}A set of vertices in $C_n$ with cardinality $m$ is \emph{maximally even} if it is of the form $J^r_{n,m}=\left\{\floor{\frac{ni+r}{m}}:i=0,1,\ldots,m-1\right\}$ for some integer $r$ with $0\leq r\leq n-1$. \end{definition}

Let us note that the $r$ in Definition \ref{definition_maximally_even} can be thought of as a rotation parameter, and reflects the fact that any rotation of a maximally even set is still maximally even. For example, $J^1_{12,5}=\{0,2,5,7,9\}$ (which is a rotation of $J^0_{12,5}$). Using Definition \ref{definition_maximally_even}, one can easily compute many instances of maximally even sets.

\section*{The Energy of a Set of Vertices}

A graph $G$ is a collection of vertices $V(G)$ together with a collection of edges $E(G)$ connecting the vertices. An edge $e$ connecting vertex $u$ to vertex $v$ is viewed as a pair $\{u,v\}$. For our purposes, we will assume all graphs are finite and simple (i.e. there will never be more than one edge between any two vertices and no edge starts and ends at the same vertex). Given two vertices $u,v\in G$ with $u\neq v$, a \emph{path from $u$ to $v$} is any finite sequence of edges $(e_1,\ldots,e_{k-1})$ such that there is a sequence of distinct vertices $(v_1,\ldots,v_k)$ with $v_1=u$, $v_k=v$, and $e_i=\{v_i,v_{i+1}\}$ for $i=1,\ldots,k-1$. The \emph{length} of a path from $u$ to $v$ is the number of edges in the path. For two vertices $u$ and $v$, the \emph{distance from $u$ to $v$} is the length of a shortest path from $u$ to $v$. We say that $G$ is \emph{connected} if for all vertices $u$ and $v$ in $G$ with $u\neq v$, there is a path from $u$ to $v$.

Let us return to one of our motivating questions: given a finite simple connected graph $G$, can we define some useful notion of ``maximally even'' set of vertices in $G$? One would hope that some characterization of maximal evenness in cycles might generalize to other graphs. However, upon a careful examination of the known characterizations of maximal evenness \cite{CloughDouthett,MR2512671,MR2408358}, one quickly realizes that all of these characterizations seem to rely on features of cycle graphs that do not hold in a broader context. Let us note that one measure of evenness, namely the sum of pairwise distances $S$ (see Definition \ref{definition_energy}), does easily generalize to finite simple connected graphs, but as is well known, this measure cannot be used to characterize maximal evenness (maximizers of the sum of distances need not be maximally even, as we will see below). 

To address this, using a simple idea from basic physics, we develop a new way to measure evenness which applies in any finite simple connected graph, and which is equivalent to maximal evenness in cycles.

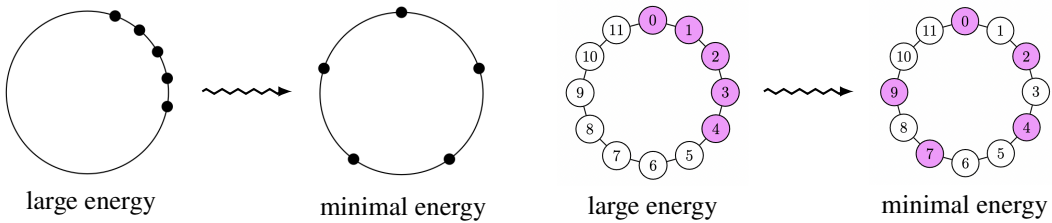
\begin{figure}[h!tbp]
\centering

\begin{tikzpicture}[x=.1\textwidth,y=.1\textwidth]
\node[inner sep=2mm] (circle_large) at (-0.5,0)
    {\includegraphics[width=1in]{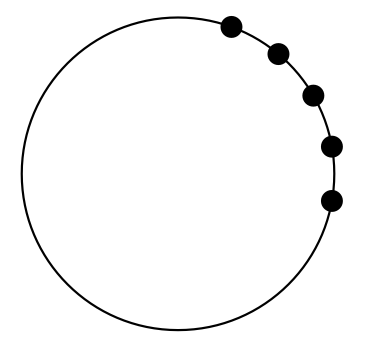}};
    \node[below = -2mm of circle_large] {\small large energy};
\node[inner sep=2mm] (circle_minimal) at (2,0)
    {\includegraphics[width=1in]{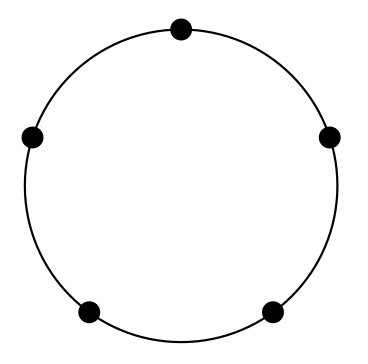}};
    \node[below = -2mm of circle_minimal] {\small minimal energy};
\draw[thick,-latex,decorate, decoration={
    zigzag,
    segment length=6,
    amplitude=.9,post=lineto,
    post length=2pt
}] (circle_large) -- (circle_minimal);

\node[inner sep=2mm] (circle_large) at (4,0)
    {\includegraphics[width=1in]{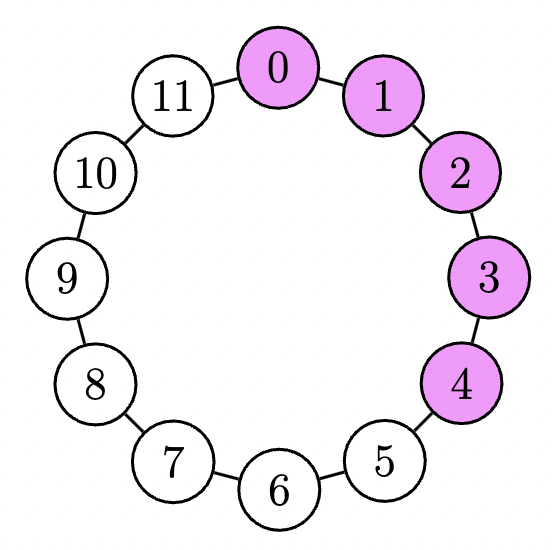}};
    \node[below = -2mm of circle_large] {\small large energy};
\node[inner sep=2mm] (circle_minimal) at (6.5,0)
    {\includegraphics[width=1in]{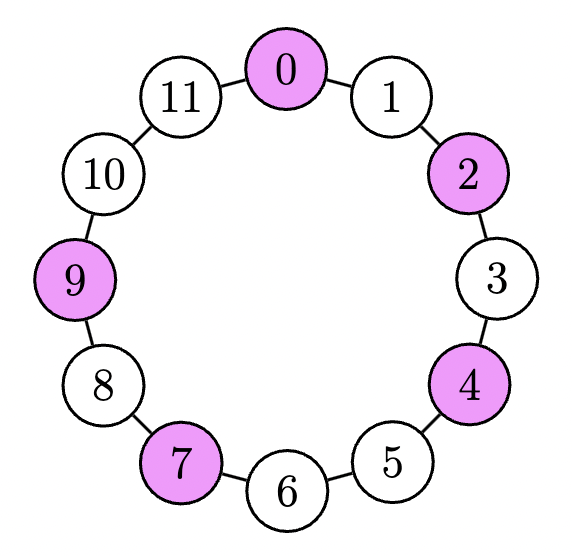}};
    \node[below = -2mm of circle_minimal] {\small minimal energy};
\draw[thick,-latex,decorate, decoration={
    zigzag,
    segment length=6,
    amplitude=.9,post=lineto,
    post length=2pt
}] (circle_large) -- (circle_minimal);
\end{tikzpicture}
\caption{Evolution of systems of unit charges on a circle and on a $12$-cycle.
}
\label{figure_equilibrium}
\end{figure}

Recall that charged particles, whose charges are of the same sign, repel one another. The electric potential energy of a system of two charged particles, each of charge $1$, i.e. massless unit point charges, is equal to $\frac{1}{d}$ where $d$ is the distance between the two particles; so when these particles are close to each other the potential energy will be large, and when they are far apart the potential energy will be small. The potential energy of a system of more than two charged particles, each of charge $1$, can be calculated by taking the sum of the reciprocals of the distances between each pair of particles in the system. So, if we place five protons on a circle in such a way that their motion is constrained to the circle, the protons will move away from each other until they reach an equilibrium configuration, as in Figure \ref{figure_equilibrium}. One way to understand this process is by using potential energy. This system of five particles on the circle starts out with some relatively large potential energy, the system evolves over time and eventually reaches a configuration of minimal potential energy. 

Our new measure of evenness can be thought of as the electrostatic potential energy of a system of unit point charges constrained to the vertices of a finite simple connected graph, using the notion of distance that comes with the graph. 


\begin{definition}\label{definition_energy}
Suppose $G$ is a finite simple connected graph and $A$ is a set of vertices in $G$ with $|A|=m$. The \emph{distance vector of $A$} is the tuple $\vec{D}(A)=(d_1,\ldots,d_\ell)$ of all distances between pairs of elements of $A$ with $d_1\leq\cdots \leq d_\ell$, where $\ell=\binom{m}{2}$. Then the \emph{energy of $A$} is $E(A)=\frac{1}{d_1}+\frac{1}{d_2}+\cdots+\frac{1}{d_\ell}$. We say that $A$ is a \emph{minimizer of $E$} (or \emph{has minimal energy}) if $E(A)$ is less than or equal to the energy of all sets of vertices in $G$ of the same cardinality as $A$. The \emph{sum of distances of $A$} is $S(A)=d_1+d_2+\cdots+d_\ell$. We say that $A$ is a \emph{maximizer of $S$} (or \emph{has maximal sum of distances}) if $S(A)$ is greater than or equal to the sum of distances of all sets of vertices in $G$ of the same cardinality as $A$. 
\end{definition}

\begin{example}[Cycles and other graphs]\label{example_cycles}
Let us calculate the sum of distances and energy of a particular subset of $C_{12}$. The maximally even set $J^0_{12,5}=\{0,2,4,7,9\}$ (see Figure \ref{figure_equilibrium}) has distance vector $\vec{D}(J^0_{12,5})=(2,2,2,3,3,4,5,5,5,5)$. This means that $S(J^0_{12,5})=36$ and $E(J^0_{12,5})=3.21\overline{6}$. Indeed, $J^0_{12,5}$ is a minimizer of $E$ and a maximizer of $S$. Let us point out that every minimizer of $E$ of size $5$ on $C_{12}$ is a rotation of $J^0_{12,5}$, while there are maximizers of $S$ of size $5$ on $C_{12}$ that are not rotations of $J^0_{12,5}$.

Similar calculations can be carried out in any finite simple connected graph using any standard mathematical software. For example, using the open source mathematical programming language SageMath, we can compute distance vectors and eneriges of all sets of vertices of size $7$ in a hypercube, and then select a minimizer as in Figure \ref{figure_hypercube_7}. Let us note that this technique indeed produces sets of vertices which are ``spaced out as evenly as possible.'' 
\end{example}

\begin{example}[Musical connections involving $C_{12}^2$]
The graph $C_{12}^2$ (see Figure \ref{figure_cyclesquared}) can be obtained from the $12$-cycle $C_{12}$ by adding an edge between any pair of vertices of $C_{12}$ which are at distance $2$ from each other in $C_{12}$. Thus, although $C_{12}$ and $C_{12}^2$ have the same vertices, distance calculations made in $C_{12}^2$ will differ from those made in $C_{12}$ because of the extra edges. This means that the two graphs $C_{12}$ and $C_{12}^2$ will have different minimizers of $E$; in some sense, the extra edges in $C_{12}^2$ lead to minimal energy sets which are ``smeared'' versions of maximally even subsets of $C_{12}$. There are exactly three minimizers of $E$ with size $7$ on $C_{12}^2$ (up to reflections and rotations), all of which are shown in Figure \ref{cyclesquared_a}-\ref{cyclesquared_c}, and indeed, none of these sets are maximally even. Hence by Theorem \ref{theorem_min_energy_and_max_evenness}, none of these sets are minimizers of $E$ on $C_{12}$. However, these sets do have musical significance. In all of the following examples, we will view sets of vertices as scales which begin at vertex $0$ and proceed in a clockwise direction unless otherwise noted. The set in Figure \ref{cyclesquared_a}, corresponds to the harmonic minor scale [C, D, E$\flat$, F, G, A$\flat$, B, C] whose sequence of steps is [W, H, W, W, H, WH, H], where W means whole step, H means half step and WH means three half steps; whereas the same set read starting on vertex $7$ and going counter clockwise corresponds to the harmonic major scale. The set in Figure \ref{cyclesquared_b} corresponds to the Chitrambari scale of Carnatic music, which has sequence of steps [W, H, W, H, W, WH, H]. The set in Figure \ref{cyclesquared_c} corresponds to the double harmonic scale in western classical music, in which the sequence of steps is [H, WH, H, W, H, WH, H]. This same scale also corresponds to indian ragas, namely the Mayamalavagowla Raga of Carnatic music and the Bhairav Raga of Hindustani classical music.


\begin{figure}[h!tbp]
\centering
\begin{minipage}[b]{0.2\textwidth} 
	\includegraphics[width=\textwidth]{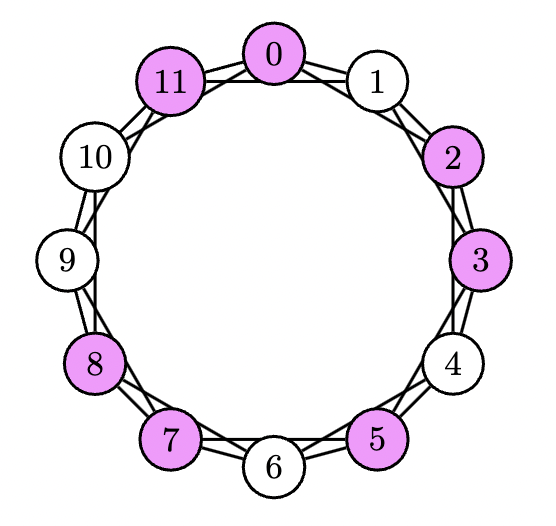}\\
    \centering\includegraphics[width=0.8\textwidth]{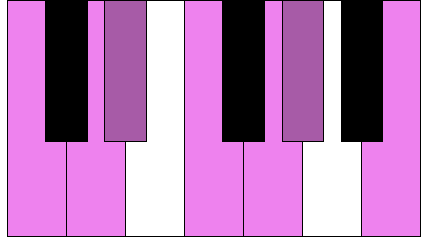}
        	\subcaption{} 
        	\label{cyclesquared_a}
\end{minipage}\quad
\begin{minipage}[b]{0.2\textwidth} 
	\includegraphics[width=\textwidth]{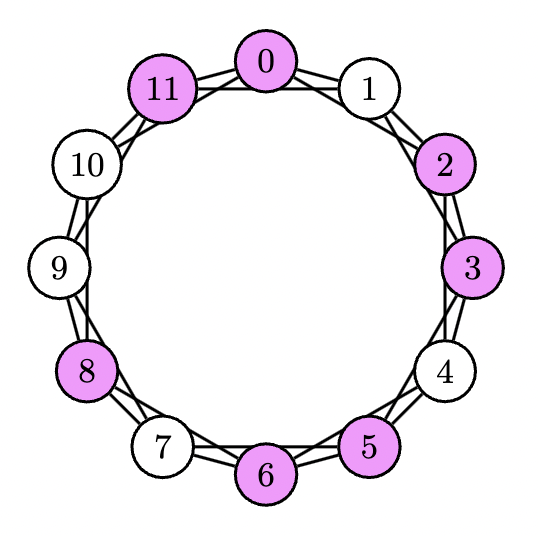}\\
    \centering\includegraphics[width=0.8\textwidth]{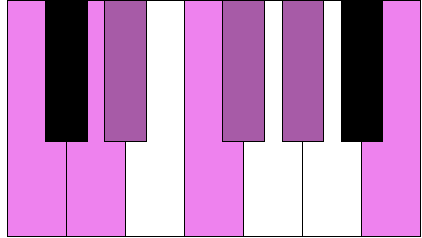}
        	\subcaption{} 
        	\label{cyclesquared_b}
\end{minipage}\quad
\begin{minipage}[b]{0.2\textwidth} 
	\includegraphics[width=\textwidth]{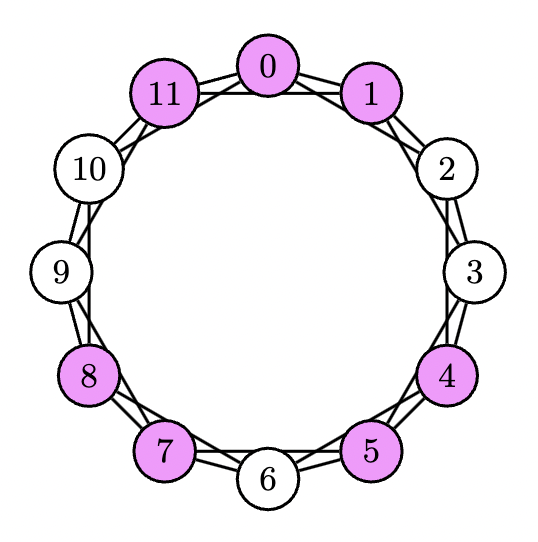}\\
    \centering\includegraphics[width=0.8\textwidth]{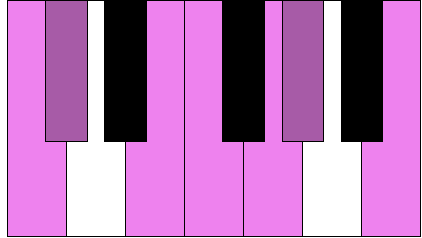}
        	\subcaption{} 
        	\label{cyclesquared_c}
\end{minipage}\quad
\begin{minipage}[b]{0.2\textwidth} 
	\includegraphics[width=\textwidth]{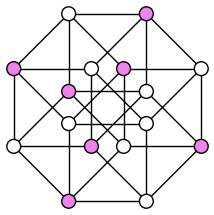}
            \subcaption{} 
        	\label{figure_hypercube_7}
\end{minipage}

~ 
\caption{The minimizers of $E$ of size $7$ on $C_{12}^2$ and the corresponding musical scales are depicted in (a)-(c), and the minimizer of $E$ of size $7$ on the hypercube appears in (d).}
\label{figure_cycle_squared}
\end{figure}

\end{example}

\begin{example}[Musical connections involving M\"obius ladder graphs: M\"obius Rhythms] 
The M\"obius ladder graph $M_{2n}$ on $2n$ vertices can be obtained from the standard ladder graph on $2n$ vertices by adding two edges joining the two ends of the ladder with a twist (see the first graph in Figure \ref{figure_mobiusladder}). The same graph $M_{2n}$ can be obtained from a cycle graph $C_{2n}$ by adding an edge between any pair of vertices $u$ and $v$ with $d(u,v)$ equal to $n$ (see $M_{16}$ in Figure \ref{figure_graphs}c). Thus, vertices which are on opposite sides of the cycle $C_{2n}$ at distance $n$ from each other will be distance $1$ from each other in $M_{2n}$. This means that minimizers of $E$ on $M_{2n}$ will tend to not contain pairs of vertices which are opposite each other whenever possible. Let us describe an interesting connection between M\"obius ladder graphs and one of the oldest Puerto Rican musical traditions of African origin called bomba. Two of the most popular rhythmic variations of this tradition are bomba yub\'a and bomba sic\'a, which are typically played on the barril de bomba (see Figure \ref{figure_barril}). The set depicted in Figure \ref{wheel_6_4} is the unique (up to reflections and rotations) minimizer of $E$ of size $4$ on $M_6$ and, when read starting on vertex $0$ and proceeding clockwise, corresponds to the bomba yuba, a 6/8-based beat which pairs well with many rhythms from Cuba and Africa \cite{Martinez}. The set shown in Figure \ref{wheel_8_4} is the unique (up to reflections and rotations) minimizer of $E$ of size $4$ on $M_8$ and corresponds to the bomba sic\'a \cite{Martinez}. The set in Figure \ref{wheel_8_6}, is one of two (up to reflections and rotations) minimizers of $E$ of size $6$ on $M_8$ and corresponds to a timbal bell (cencerro) pattern often played in bomba \cite{Mauleon}.

\begin{figure}[h!tbp]
\centering
\begin{minipage}[b]{0.18\textwidth}
\centering
	\includegraphics[width=0.8\textwidth]{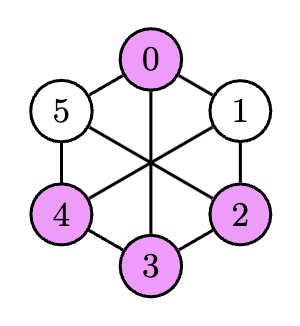}
        	\subcaption{} 
        	\label{wheel_6_4}
\end{minipage}\quad
\begin{minipage}[b]{0.18\textwidth} 
	\includegraphics[width=\textwidth]{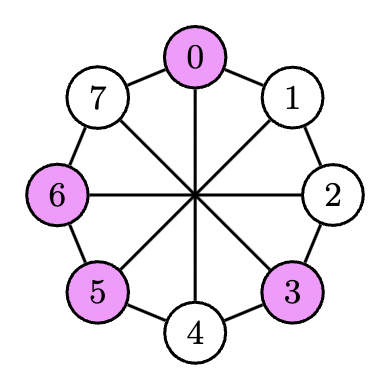}
        	\subcaption{} 
        	\label{wheel_8_4}
\end{minipage}\quad
\begin{minipage}[b]{0.18\textwidth} 
	\includegraphics[width=\textwidth]{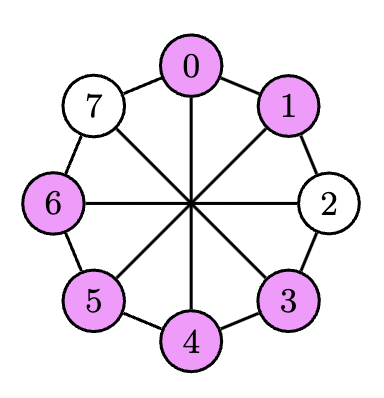}
        	\subcaption{} 
        	\label{wheel_8_6}
\end{minipage}\quad
\begin{minipage}[b]{0.15\textwidth} 
	\includegraphics[width=\textwidth]{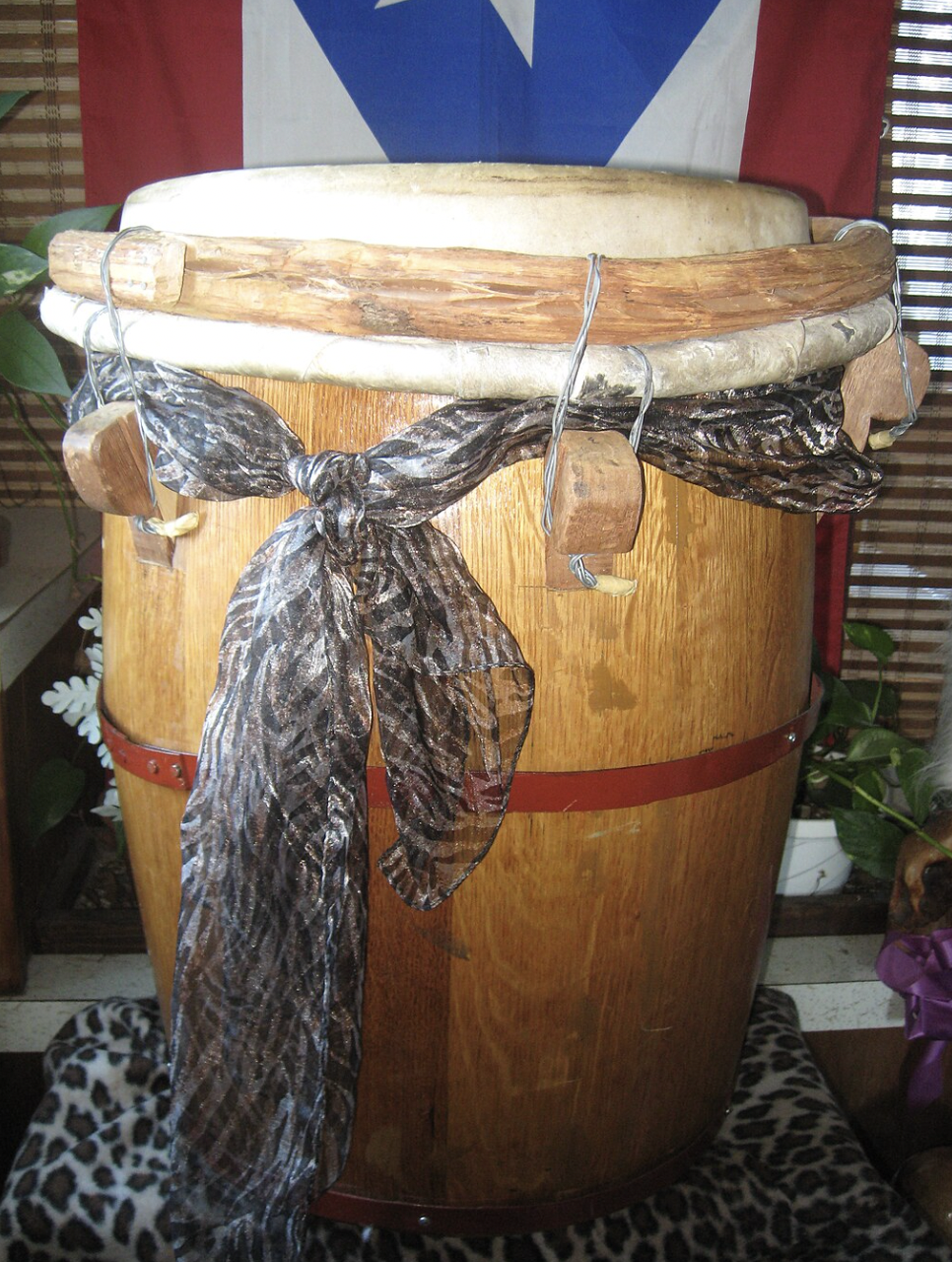}
        	\subcaption{\cite{Yolydia}} 
        	\label{figure_barril}
\end{minipage}
~ 
\caption{Some minimizers of $E$ on M\"obius ladders and a barril de bomba.}
\label{figure_wheel}
\end{figure}

\end{example}

\section*{Minimal Energy, Maximal Evenness, and Majorization}

The examples discussed above suggest that the notion of minimal energy is a good candidate for a possible generalization of maximal evenness from $n$-cycles to arbitrary finite simple connected graphs. In this section we present more evidence along these lines in the form of several general theorems, proofs of which are essentially combinatorial and will appear in a forthcoming article. Along the way we introduce another previously unknown technique for measuring the evenness of a set of vertices of a graph: \emph{majorization}. The first theorem we would like to present (which is similar to but distinct from \cite[Theorem 7]{MR2366388}) says that when we restrict our attention to $n$-cycles, minimal energy is equivalent to maximal evenness. Thus, our notion of minimal energy set of vertices in a graph is a direct generalization of Clough and Douthett's maximal evenness.

\begin{theorem}\label{theorem_min_energy_and_max_evenness}
For any positive integer $n$, a set of vertices $A$ of an $n$-cycle $C_n$ has minimal energy if and only if it is maximally even.
\end{theorem}

It is well known that the complement of a maximally even set is maximally even (see \cite[Theorem 3.3]{CloughDouthett}, \cite[Corollary 1]{MR2512671}, or \cite[Theorem 3.1]{MR2537012}). Since we are claiming that the notion of minimal energy set provides a good generalization of the notion of maximal evenness, it is desirable that the complement of a minimal energy set also has minimal energy. However, this is not true in general. For example, the subset $A$ of the path graph $P_4$ on four vertices consisting of the two endpoints has minimal energy $E(A)=\frac{1}{3}$, whereas its complement $A^c=V(P_4)\setminus A$ has energy $E(A^c)=1$, which is clearly not minimal. This being said, there is a wide class of graphs, namely the distance degree regular graphs introduced by Bloom et al. \cite{MR0634519}, in which complements of minimal energy sets always have minimal energy. Given a vertex $v$ in $G$, the \emph{distance vector of $v$} is the vector $\vec{D}(v)=(d_1,\ldots,d_L)$ of length $L=|V(G)|-1$ of all distances from $v$ to other vertices in the graph, where $d_1\leq \cdots\leq d_L$. A graph $G$ is \emph{distance degree regular} if whenever $u,v\in G$ we have $\vec{D}(u)=\vec{D}(v)$. Notice that squares of cycles, M\"obius ladders, and the hypercube are all distance degree regular graphs, so by the following theorem the complements of the sets in Figure \ref{figure_cycle_squared}, and Figure \ref{figure_wheel} are also minimizers of $E$.


\begin{theorem}\label{theorem_complements}
Suppose $G$ is a finite connected simple graph which is distance degree regular. Then a set of vertices $A$ of $G$ has minimal energy if and only if its complement $V(G)\setminus A$ has minimal energy. 
\end{theorem}

One of the main tools used in the proof of Theorem \ref{theorem_min_energy_and_max_evenness} is the idea of \emph{majorization} (see Definition \ref{definition_majorization}), which provides us with another way to measure the evenness of a set of vertices in finite connected graphs. Since its origin over one-hundred years ago in Muirhead's 1902 study of mathematical inequalities \cite{muirhead1902} and in Lorenz's 1905 study of wealth inequality \cite{lorenz1905}, majorization has played an important role in many areas of research \cite{MR2759813}, such as linear algebra \cite{MR0995373}, operator theory \cite{MR0203464}, and brownian motion \cite{MR0474525}. Let us consider an example. Suppose we need to distribute twenty apples among four people. We'd like to have a way of comparing the equity of such distributions. For example, one such distribution can be denoted by $x=(1,4,5,10)$ and another by $y=(1,3,6,10)$, where the $i$th coordinate of each tuple is the number of apples received by person $i$. Which distribution is more equitable? The idea is that $x$ is more equitable than $y$ because $x$ can be obtained from $y$ via a \emph{Robinhood transfer}, meaning that we can start with $y=(1,3,6,10)$ and transfer one apple from some person to some other person who has fewer apples, thus obtaining $x=(1,3+1,6-1,10)=(1,4,5,10)$. The tuple $(1,4,5,10)$ is \emph{majorized} by (i.e. is more equitable than) $(1,3,6,10)$. Of course, the most equitable distribution in this example would be $(5,5,5,5)$, which can be obtained by performing a finite sequence of Robinhood transfers on any vector $(y_1,y_2,y_3,y_4)$ of nonnegative integers whose entries sum to twenty.

The general definition of majorization, which we discuss now, provides an efficient way to compare the ``equity'' of $\ell$-tuples when viewed as wealth distributions. Given an $\ell$-tuple $x=(x_1,\ldots,x_\ell)$ we let $x^\downarrow=(x_1^\downarrow,\ldots,x_\ell^\downarrow)$ denote the $\ell$-tuple which has the same entries as $x$ but $x_1^\downarrow\geq\cdots\geq x_\ell^\downarrow$. 

\begin{definition}\label{definition_majorization}
Suppose $x=(x_1,\ldots,x_\ell)$ and $y=(y_1,\ldots,y_\ell)$ are $\ell$-tuples of real numbers. We say that $x$ is \emph{majorized} by $y$ and write $x\prec y$ if and only if (1) for $k=1,\ldots,\ell-1$ we have $\sum_{i=1}^kx_i^\da\leq\sum_{i=1}^ky_i^\da$ and (2) $\sum_{i=1}^\ell x_i^\da=\sum_{i=1}^\ell y_i^\da$.
\end{definition}

But what does majorization have to do with the energy of a set of vertices? How is majorization used to prove Theorem \ref{theorem_min_energy_and_max_evenness}? The following theorem, due to the great Issai Schur \cite{schur1923} and independently to the eminent team of mathematicians consisting of Hardy, Littlewood, and P\'olya \cite{HLP}, provides a connection between majorization and the energy of subsets of graphs, and in particular cycles.

\begin{theorem}[{Schur \cite{schur1923}; Hardy, Littlewood, and Polya \cite{HLP}}]\label{theorem_HLP}
Suppose $I\subseteq\R$ is an interval, and $\phi:I^\ell\to\R$ is defined by $\phi(x)=\sum_{i=1}^\ell \frac{1}{x_i}.$ Then $x\prec y$ implies $\phi(x)\leq \phi(y)$. Moreover, if $x\prec y$ and $x$ is not a permutation of $y$ then it follows that $\phi(x)<\phi(y)$.
\end{theorem}

The intuition provided by Theorem \ref{theorem_HLP} is that if $A$ and $B$ are sets of vertices in a finite connected graph $G$, and if $\vec{D}(A)$ (as in Definition \ref{definition_energy}) is majorized by $\vec{D}(B)$, then $\vec{D}(A)$ is more equitable than $\vec{D}(B)$, and that the vertices in $A$ are spaced out ``more evenly'' than those in $B$. This intuition is further justified by the following.

\begin{theorem}\label{theorem_majorization}
Suppose $n\geq 3$ is an integer and $1\leq m\leq n$. Let $\mathcal{S}_{n,m}$ be the collection of all sets of vertices of $C_n$ of cardinality $m$ which maximize the sum of distances $S$. Then a set of vertices $A\subseteq V(C_n)$ of size $m$ is maximally even if and only if $\vec{D}(A)\prec \vec{D}(B)$ for all $B\in \mathcal{S}_{n,m}$.
\end{theorem}

The sum of distances, the notion of energy, and also majorization all provide different ways to measure the evenness of a set of vertices. Rather than declaring one measure of evenness preferable to another, we tend to view these three notions as working well together, with part of the picture remaining hidden if one or more of the measures is discarded. Let us consider a few examples that involve musical connections and illustrate the way in which the nonlinear structure of the majorization order fits together with the sum of distances and the notion of energy.



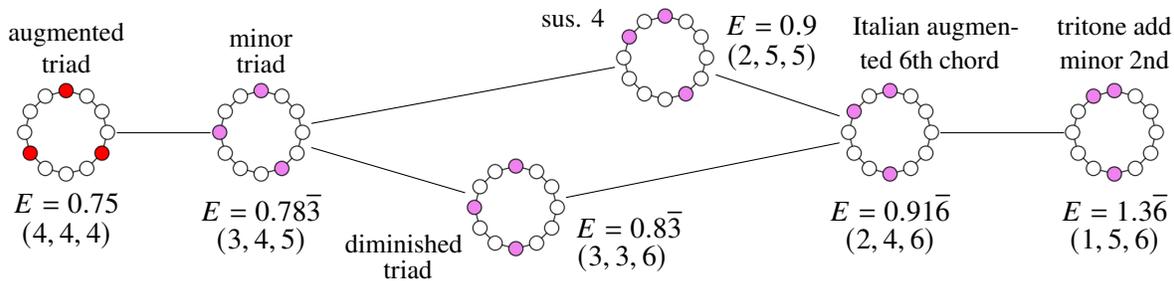
\begin{figure}[h!tbp]
	\centering
\begin{tikzpicture}[x=.12\textwidth,y=.1\textwidth]
\node[inner sep=0pt] (wi_12_3_0) at (0,0)
    {\includegraphics[width=.08\textwidth]{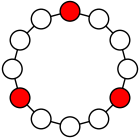}};
    \node[above = 0.9em of wi_12_3_0,align=center] {\small augmented};
    \node[above = 0mm of wi_12_3_0] {\small triad};
    \node[below = 0mm of wi_12_3_0] {$E=0.75$};
    \node[below = 1em of wi_12_3_0] {$(4,4,4)$};
\node[inner sep=0pt] (wi_12_3_1) at (1.3,0)
    {\includegraphics[width=.08\textwidth]{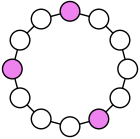}};
    \node[above = 0.9em of wi_12_3_1] {\small minor};
    \node[above = 0mm of wi_12_3_1] {\small triad};
    \node[below = 0mm of wi_12_3_1] {$E=0.78\overline{3}$};
    \node[below = 1.2em of wi_12_3_1] {$(3,4,5)$};
\node[inner sep=0pt] (wi_12_3_2) at (3,-.06\textwidth)
    {\includegraphics[width=.08\textwidth]{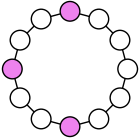}};
    \node[below left = -4.5mm and -1mm of wi_12_3_2] {\small diminished};
    \node[below left = -1mm and 3mm of wi_12_3_2] {\small triad};
    \node[below right = -2em and 0mm of wi_12_3_2] {$E=0.8\overline{3}$};
    \node[below right = -0.8em and 0mm of wi_12_3_2] {$(3,3,6)$};
\node[inner sep=0pt] (wi_12_3_3) at (4,.06\textwidth)
    {\includegraphics[width=.08\textwidth]{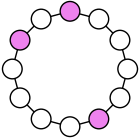}};
    \node[above left = -1em and 0 em of wi_12_3_3] {\small sus. 4};
    \node[above right = -1.3em and 0mm of wi_12_3_3] {$E=0.9$};
    \node[above right = -2.5em and 0mm of wi_12_3_3] {$(2,5,5)$};
\node[inner sep=0pt] (wi_12_3_4) at (5.5,0)
    {\includegraphics[width=.08\textwidth]{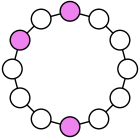}};
    \node[above right = 4.25mm and -13mm of wi_12_3_4] {\small Italian augmen-};
    \node[above right = 0.5mm and -12mm of wi_12_3_4] {\small ted 6th chord};
    \node[below = 0mm of wi_12_3_4] {$E=0.91\overline{6}$};
    \node[below = 1.2em of wi_12_3_4] {$(2,4,6)$};
\node[inner sep=0pt] (wi_12_3_5) at (7,0)
    {\includegraphics[width=.08\textwidth]{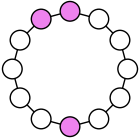}};
    \node[above = 5mm of wi_12_3_5] {\small tritone add};
    \node[above = 0.5mm of wi_12_3_5] {\small minor 2nd};
    \node[below = 0mm of wi_12_3_5] {$E=1.3\overline{6}$};
    \node[below = 1.2em of wi_12_3_5] {$(1,5,6)$};
\draw[-] (wi_12_3_0) -- (wi_12_3_1);
\draw[-] (wi_12_3_1) -- (wi_12_3_2);
\draw[-] (wi_12_3_1) -- (wi_12_3_3);
\draw[-] (wi_12_3_2) -- (wi_12_3_4);
\draw[-] (wi_12_3_3) -- (wi_12_3_4);
\draw[-] (wi_12_3_4) -- (wi_12_3_5);
    
\end{tikzpicture}

\caption{Majorization order on distance vectors of three-note chords that maximize the sum of distances.}
	\label{figure_majorization_12_3}
\end{figure}

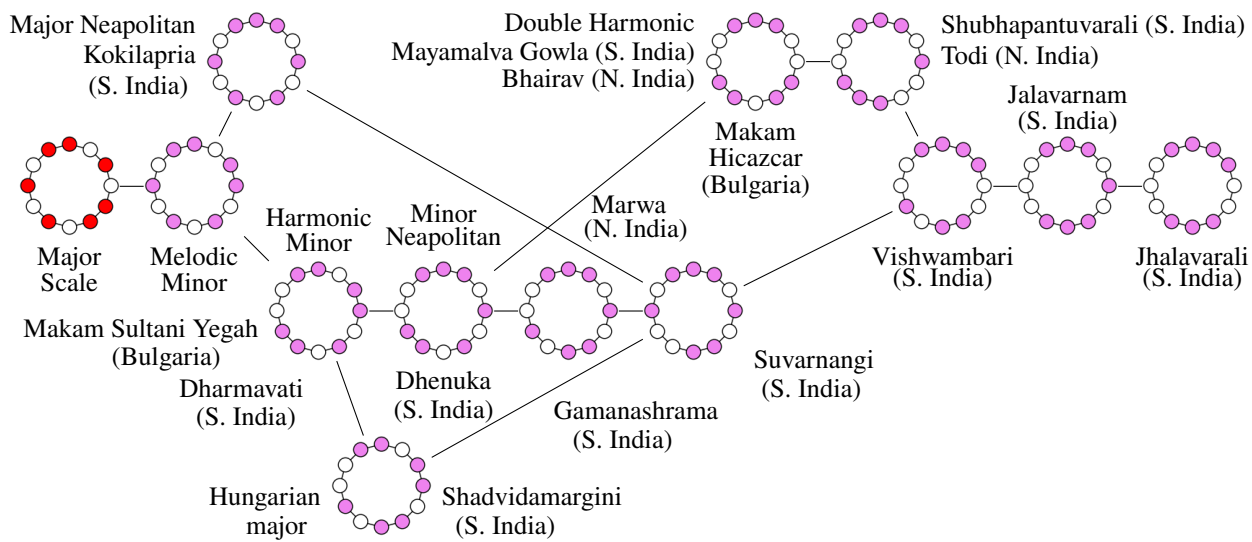
\begin{figure}[h!tbp]
	\centering
\begin{tikzpicture}[x=.1\textwidth,y=.1\textwidth]
\node[inner sep=0pt] (wi_12_7_0) at (0,0)
    {\includegraphics[width=.08\textwidth]{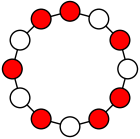}};
    \node[below = 0mm of wi_12_7_0] {\small Major};
    \node[below = 3.5mm of wi_12_7_0] {\small Scale};
\node[inner sep=0pt] (wi_12_7_1) at (1,0)
    {\includegraphics[width=.08\textwidth]{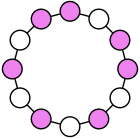}};
    \node[below = 0mm of wi_12_7_1] {\small Melodic};
    \node[below = 0.9em of wi_12_7_1] {\small Minor};
\node[inner sep=0pt] (wi_12_7_2) at (1.5,1)
    {\includegraphics[width=.08\textwidth]{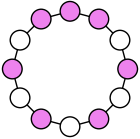}};
    \node[above left = -5mm and 0mm of wi_12_7_2] {\small Major Neapolitan};
    \node[above left = -9mm and 0mm of wi_12_7_2] {\small Kokilapria};
    \node[above left = -13mm and 1mm of wi_12_7_2] {\small (S. India)};
\node[inner sep=0pt] (wi_12_7_3) at (2,-1)
    {\includegraphics[width=.08\textwidth]{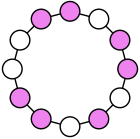}};
    \node[above = 3.3mm of wi_12_7_3] {\small Harmonic};
    \node[above = 0mm of wi_12_7_3] {\small Minor};
    \node[below left = -7mm and 0mm of wi_12_7_3] {\small Makam Sultani Yegah};
    \node[below left = -3.2mm and 5mm of wi_12_7_3] {\small (Bulgaria)};
    \node[below left = 1mm and -6mm of wi_12_7_3] {\small Dharmavati};
    \node[below left = 4.8mm and -5mm of wi_12_7_3] {\small (S. India)};
\node[inner sep=0pt] (wi_12_7_4) at (2.5,-2.4)
    {\includegraphics[width=.08\textwidth]{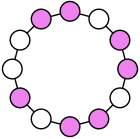}};
    \node[below left = -7.8mm and 0mm of wi_12_7_4] {\small Hungarian};
    \node[below left = -4mm and 1.8mm of wi_12_7_4] {\small major};
    \node[below right = -7.8mm and 0mm of wi_12_7_4] {\small Shadvidamargini};
    \node[below right = -4mm and 1.8mm of wi_12_7_4] {\small (S. India)};
\node[inner sep=0pt] (wi_12_7_5) at (3,-1)
    {\includegraphics[width=.08\textwidth]{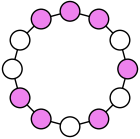}};
    \node[above = 3.8mm of wi_12_7_5] {\small Minor};
    \node[above = 0mm of wi_12_7_5] {\small Neapolitan};
    \node[below = 0mm of wi_12_7_5] {\small Dhenuka};
    \node[below = 3.8mm of wi_12_7_5] {\small (S. India)};
\node[inner sep=0pt] (wi_12_7_6) at (4,-1)
    {\includegraphics[width=.08\textwidth]{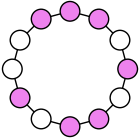}};
    \node[below right = 4mm and -10mm of wi_12_7_6] {\small Gamanashrama};
    \node[below right = 7.8mm and -7.5mm of wi_12_7_6] {\small (S. India)};
    \node[above right = 5mm and -5mm of wi_12_7_6] {\small Marwa};
    \node[above right = 1.2mm and -6mm of wi_12_7_6] {\small (N. India)};
\node[inner sep=0pt] (wi_12_7_7) at (5,-1)
    {\includegraphics[width=.08\textwidth]{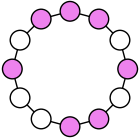}};
    \node[below right = -2.8mm and 0mm of wi_12_7_7] {\small Suvarnangi};
    \node[below right = 1mm and 1mm of wi_12_7_7] {\small (S. India)};
\node[inner sep=0pt] (wi_12_7_8) at (5.5,1)
    {\includegraphics[width=.08\textwidth]{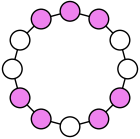}};
    \node[above left = -4mm and 0mm of wi_12_7_8] {\small Double Harmonic};
    \node[above left = -8.5mm and 0mm of wi_12_7_8] {\small Mayamalva Gowla (S. India)};
    \node[above left = -12mm and 0mm of wi_12_7_8] {\small Bhairav (N. India)};
    \node[below = 0mm of wi_12_7_8] {\small Makam};
    \node[below = 3.5mm of wi_12_7_8] {\small Hicazcar};  
    \node[below = 7mm of wi_12_7_8] {\small (Bulgaria)};  
\node[inner sep=0pt] (wi_12_7_9) at (6.5,1)
    {\includegraphics[width=.08\textwidth]{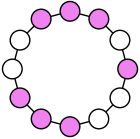}};
    \node[above right = -5mm and 0mm of wi_12_7_9] {\small Shubhapantuvarali (S. India)};
    \node[above right = -9mm and 0mm of wi_12_7_9] {\small Todi (N. India)};
\node[inner sep=0pt] (wi_12_7_10) at (7,0)
    {\includegraphics[width=.08\textwidth]{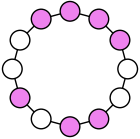}};
    \node[below = 0mm of wi_12_7_10] {\small Vishwambari};
    \node[below = 2.8mm of wi_12_7_10] {\small (S. India)};
\node[inner sep=0pt] (wi_12_7_11) at (8,0)
    {\includegraphics[width=.08\textwidth]{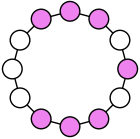}};
    \node[above = 3mm of wi_12_7_11] {\small Jalavarnam};
    \node[above = -1mm of wi_12_7_11] {\small (S. India)};
\node[inner sep=0pt] (wi_12_7_12) at (9,0)
    {\includegraphics[width=.08\textwidth]{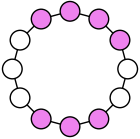}};
    \node[below = 0mm of wi_12_7_12] {\small Jhalavarali};
    \node[below = 2.8mm of wi_12_7_12] {\small (S. India)};

\draw[-] (wi_12_7_0) -- (wi_12_7_1);
\draw[-] (wi_12_7_1) -- (wi_12_7_2);
\draw[-] (wi_12_7_1) -- (wi_12_7_3);
\draw[-] (wi_12_7_3) -- (wi_12_7_4);
\draw[-] (wi_12_7_3) -- (wi_12_7_5);
\draw[-] (wi_12_7_5) -- (wi_12_7_6);
\draw[-] (wi_12_7_6) -- (wi_12_7_7);
\draw[-] (wi_12_7_2) -- (wi_12_7_7);
\draw[-] (wi_12_7_4) -- (wi_12_7_7);
\draw[-] (wi_12_7_5) -- (wi_12_7_8);
\draw[-] (wi_12_7_8) -- (wi_12_7_9);
\draw[-] (wi_12_7_9) -- (wi_12_7_10);
\draw[-] (wi_12_7_7) -- (wi_12_7_10);
\draw[-] (wi_12_7_10) -- (wi_12_7_11);
\draw[-] (wi_12_7_11) -- (wi_12_7_12);

\end{tikzpicture}

\caption{Majorization order on distance vectors of seven-note scales that maximize the sum of distances.}
	\label{figure_majorization_12_7}
\end{figure}

\begin{example} Using the notation of Theorem \ref{theorem_majorization}, each set in $\mathcal{S}_{12,3}$ has cardinality $3$, has a distance vector of length $\binom{3}{2}=3$ and corresponds to a musical chord (see Figure \ref{figure_majorization_12_3}), whereas each set in $\mathcal{S}_{12,7}$ has cardinality $7$, has a distance vector of length $\binom{7}{2}=21$ and corresponds to a musical scale (see Figure \ref{figure_majorization_12_7}). Let us consider the majorization order on the distance vectors of sets in $\mathcal{S}_{12,3}$ and $\mathcal{S}_{12,7}$. 

In both Figure \ref{figure_majorization_12_3} and Figure \ref{figure_majorization_12_7}, the set on the left is maximally even, has minimal energy and its distance vector is majorized by the distance vector of every other set in the diagram. The energies of the displayed sets increase from left to right and each line indicates a majorization between the distance vectors of the displayed sets, where one can ``compose'' majorization by transitivity. For example in Figure \ref{figure_majorization_12_3}, $(3,4,5)\prec(2,5,5)$ and $(3,4,5)\prec(2,4,6)$. It is worth noting that the majorization order on both sets of distance vectors is nonlinear. For example, in Figure \ref{figure_majorization_12_3}, the distance vector $(3,3,6)$ is not comparable to $(2,5,5)$, i.e. $(3,3,6)\not\prec (2,5,5)$ and $(2,5,5)\not\prec(3,3,6)$. Let us also point out that reading Figure \ref{figure_majorization_12_7} from left to right, each time a new minor third (i.e. an interval consisting of three half steps) appears, a split in the majorization order occurs.



\end{example}


\section*{Summary and Conclusions}

The notion of energy of a set of vertices serves as a suitable measure of evenness in any finite simple connected graph. Indeed, due to the mathematical theorems discussed above, the notion of minimal energy set of vertices is a faithful mathematical generalization of maximal evenness. We see the connections between M\"{o}bius ladder graphs and bomba, between $C_{12}^2$ and various scales as well as between the majorization order and scales, as providing evidence that our notion of energy, the sum of distances and majorization are worth studying in both the cyclic and the non-cyclic more geometrically complex context.


\section*{Acknowledgements}

Special thanks to Jason Laferrera for creating software that allowed us to listen to minimal energy sets.

    
{\setlength{\baselineskip}{13pt} 
\raggedright				
\bibliographystyle{bridges}
\bibliography{bridgesbib}
} 
   
\end{document}